\documentclass[12pt]{amsart}
\usepackage{amscd, amssymb,latexsym,amsmath, amscd, amsmath}
\usepackage[all]{xy}

\newtheorem{theorem}{Theorem}[section]
\newtheorem{lemma}[theorem]{Lemma}
\newtheorem{prop}[theorem]{Proposition}
\newtheorem{remark}[theorem]{Remark}
\newtheorem{thm}[theorem]{Theorem}
\newtheorem{cor}[theorem]{Corollary}

\begin{document}

\title[On the ${\rm SL(2)}$ period integral]
{On the $SL(2)$ period integral}
\author{U. K. Anandavardhanan}
\address{Department of Mathematics, Indian Institute of Technology Bombay, Bombay - 400 076, INDIA.}
\email{anand@math.iitb.ac.in}

\author{Dipendra Prasad}
\address{Tata Institute of Fundamental 
Research, Homi Bhabha Road, Bombay - 400 005, INDIA.}
\email{dprasad@math.tifr.res.in}


\date{}

\begin{abstract}
Let $E/F$ be a quadratic extension of number fields. For a cuspidal
representation $\pi$ of ${\rm SL}_2({\Bbb A}_E)$, we study in this paper 
the integral of functions in $\pi$ on ${\rm SL}_2(F)\backslash {\rm SL}_2({\Bbb A}_F)$. 
We characterize the non-vanishing of these
integrals, called period integrals,
 in terms of $\pi$ having a  Whittaker model with respect to characters
of $E\backslash {\Bbb A}_E$ which are trivial on ${\Bbb A}_F$.
We show that the period 
integral in general is not a product of local invariant functionals,
and find a necessary and sufficient condition when it is.
We exhibit cuspidal representations of ${\rm SL}_2({\Bbb A}_E)$ whose period
integral vanishes identically while each local
constituent admits an ${\rm SL}_2$-invariant linear functional. Finally, 
 we construct an automorphic representation $\pi$
on ${\rm SL}_2({\Bbb A}_E)$ which is
 abstractly ${\rm SL}_2({\Bbb A}_F)$ distinguished
but for which none of the elements in the global $L$-packet 
determined by it  is  distinguished by ${\rm SL}_2({\Bbb A}_F)$.
\end{abstract}

\maketitle

\setcounter{tocdepth}{1}

\tableofcontents

\section{Introduction}

Let $F$ be a number field and ${\Bbb A}_F$ its ad\`ele ring. Let $G$ be
a reductive algebraic group over $F$ and $H$ a reductive subgroup of $G$
over $F$. 
Assume that the center $Z_H$ of $H$ is contained in the center $Z_G$ of $G$,
a condition that holds in the cases we study in this paper.
For an automorphic form $\phi$ on $G({\Bbb A}_F)$ on which $Z_H({\Bbb A}_F)$
acts trivially, 
the {\it period integral}
of $\phi$ with respect to $H$ is defined to be the integral
(when convergent, which is the case  if $\phi$ is cuspidal)
$${\mathcal P}(\phi)=\int_{H(F)Z_H({\Bbb A}_F)
\backslash H({\Bbb A}_F)}\phi(h)dh,$$
where $dh$ is the natural measure on $H(F)Z_H({\Bbb A}_F)
\backslash H({\Bbb A}_F)$,
the so called Tamagawa measure. An automorphic representation $\pi$
of $G({\Bbb A}_F)$ is said to be globally distinguished with respect
to $H$ if this period integral
is nonzero for some $\phi \in \pi$. More
generally, if $\chi$ is a one-dimensional representation of
$H({\Bbb A}_F)$ trivial on $H(F)$ such that 
$Z_H({\Bbb A}_F)$ acts trivially on
$\phi(h)\chi(h)$, and 
$$\int_{H(F)Z_H({\Bbb A}_F)\backslash H({\Bbb A}_F)}\phi(h)\chi(h)dh,$$
is nonzero for some $\phi \in \pi$, then $\pi$ is said to be 
$\chi$-distinguished.

The notion of distinguishedness has been extensively studied (and is
being studied) both locally and globally  for $G = {\rm GL}_n({\Bbb A}_E)$, 
and $H = {\rm GL}_n({\Bbb A}_F)$, $E$ a quadratic extension of a number field $F$
\cite{flicker1,flicker2,hakim1,harder}. 
In this case (just as in the local case), distinguishedness implies that 
$$\pi^{\sigma} \cong \pi^{\vee},$$
where $\sigma$ is the nontrivial automorphism of $E/F$, and 
is equivalent to the Asai $L$-function
of $\pi$, denoted by say $L(s, r(\pi))$, having a pole at $s =1$. From the factorization
of $L$-functions:
$$L(s,\pi \times \pi^\sigma)=L(s,r(\pi))L(s,r(\pi)\otimes \omega_{E/F}),$$
it follows that if
$\pi=\otimes \pi_v$, $v$ running over all places of $E$, is an
automorphic representation of $G = {\rm GL}_n({\Bbb A}_E)$, with all the
local components $\pi_v$ distinguished, then $\pi$ is either distinguished,
or $\omega_{E/F}$ distinguished, where $\omega_{E/F}$ is the
quadratic character of ${\Bbb A}_F^*$ associated to the extension
$E/F$. Thus local distinguishedness for ${\rm GL}_n$ `almost' implies
global distinguishedness. 
Observe furthermore that $\pi$ cannot be both distinguished 
and $\omega_{E/F}$-distinguished
(as it would then 
contribute a pole of order 2 to the Rankin product $L$-function 
$L(s, \pi \otimes \pi^\sigma)$).

In an earlier work, the authors had studied the distinguishedness property
for ${\rm SL}_2$ in the local case. We carry out the global analysis of this 
case here. 
Since for ${\rm GL}_n$ as mentioned earlier, locally distinguished automorphic
representations are either distinguished, or $\omega_{E/F}$ distinguished,
one is led to ask whether locally distinguished automorphic
representations of ${\rm SL}_2({\Bbb A}_E)$ are globally  distinguished
by ${\rm SL}_2({\Bbb A}_F)$. We show in this paper that this is not the
case. 

If $\pi$ is an automorphic representation of $G({\Bbb A}_F)$ for $G$ 
a general reductive group over a number field $F$, then 
$\pi$ factorizes as  $\pi=\otimes \pi_v$, $v$ running over all places of $F$. 
For an algebraic subgroup $H$ of $G$ defined over $F$, the period
integral
$\phi \mapsto {\mathcal P}(\phi)$ is
an $H({\Bbb A}_F)$-invariant form on $\pi$. If one knows that the space of 
$H(F_v)$-invariant forms on an irreducible representation
of $G(F_v)$ is at most one dimensional for any place $v$ of $F$, 
then the invariant form
$\phi \mapsto {\mathcal P}(\phi)$ is a ``product'' of local invariant 
forms times a global constant which one expects to be intimately 
connected with special values of automorphic $L$-functions 
associated to $\pi$.  

Recently, a very interesting case has been studied by Jacquet in \cite
{jacquet3} where the space of $H(F_v)$-invariant forms on an irreducible 
admissible representation of $G(F_v)$ is not always one dimensional 
but for which the functional $\phi \mapsto {\mathcal P}(\phi)$ is nevertheless
expressible as a product of local factors. Jacquet's example is for the case:
$(G,H)=({\rm Res}_{E/F}{\rm GL}_3,{\rm U}_3)$. We have recently learnt
that Jacquet has generalized this work to $GL_n$.

In the earlier work \cite{anand} where the authors analyzed the situation for 
$G={\rm Res}_{E/F}{\rm SL}_2$ and $H={\rm SL}_2$ locally, it was found that 
multiplicity one can fail for the space of $H$-invariant forms even for an 
irreducible admissible representation of $G$  which is
supercuspidal. We analyze in this paper 
whether the period integral
is factorizable in this case. We find that this is so if 
the automorphic representation
$\pi$ is not monomial,
and also in the case when it is monomial and comes from 3 quadratic
extensions of $E$ of which only one is Galois over $F$; in other
cases, the period integral is not factorizable. What is most appealing
about this result is that it is the exact global analogue (interpreted
via Galois theory) of the local results obtained in [A-P] about the dimension
of the space of $SL_2(k)$-invariant forms for a representation of 
$SL_2(K)$ where $K$ is a quadratic extension of a non-archimedean local
field $k$.  

In trying to understand representations of ${\rm GL}_2({\Bbb A}_E)$
which are distinguished by ${\rm SL}_2({\Bbb A}_F)$, we are naturally
led to investigate a related concept, which we call {\it 
pseudo-distinguishedness}. They are studied in section 7.

The main  results proved in this paper are theorems 4.2, 5.2,
6.8 and 8.2. In section 2 we take up some preliminary results
about the structure of $L$-packets for $SL_2$. In particular,
a rather simple proof is provided for the stability of non-monomial 
representations for $SL_2$ and
more generally of primitive representations for $SL_n$; see 
Lemma \ref{lemma2.4}, and Remark \ref{remark2.7}.

We end the introduction by mentioning that in this paper we have
constructed examples  of
automorphic representations
$\pi=\otimes \pi_v$  of ${\rm SL}_2({\Bbb A}_E)$ which are
 abstractly ${\rm SL}_2({\Bbb A}_F)$-distinguished but which are
not globally  distinguished,
and also constructed examples  of
automorphic representations
$\pi=\otimes \pi_v$  of ${\rm SL}_2({\Bbb A}_E)$ which are
 abstractly ${\rm SL}_2({\Bbb A}_F)$-distinguished but for which
none of the elements in the global $L$-packet is globally
distinguished. Thus local distinguishedness fails to guarantee global 
distinguishedness even at the level of $L$-packets. The reader will note that
the above statements have the flavor of Blasius' results \cite{blasius}.
But we have not achieved a real understanding
of this phenomenon in this paper. Perhaps there is a certain
multiplicity formula in the spirit of Labesse-Langlands, cf.
\cite{labesse},
which determines when a member of an $L$-packet determined by $\pi$
has non-vanishing period integral on  ${\rm SL}_2(F)
\backslash {\rm SL}_2({\Bbb A}_F)$; this we have not been able to do here.

\vspace{2mm}

\noindent{\bf Acknowledgement:} The authors thank the referee for 
a careful reading and detailed comments on the paper.

\section{Some lemmas about $L$-packets on ${\rm SL}_2$}

In this section we recall for the reader's convenience some standard results 
about ${\rm SL}_2$. We supply the proofs to emphasize their elementary 
nature; see also \cite{lapid}.

\indent
For an irreducible representation
$\pi$ of a normal subgroup $N$ of a group $G$ and an element $g$ of $G$, we let
${}^g \pi$ denote the representation $n \mapsto \pi(gng^{-1})$ of $N$.

\begin{lemma}\label{lemma2.1}
If $\pi_1=\otimes \pi_{1,v}$ and $\pi_2=\otimes \pi_{2,v}$ are two 
cuspidal automorphic
representations of ${\rm SL}_2({\Bbb A}_F)$ which are in the same $L$-packet, 
i.e., $\pi_1$ and $\pi_2$ are irreducible subrepresentations of an automorphic
representation $\tilde{\pi}$ of ${\rm GL}_2({\Bbb A}_F)$,
then there exists $g \in {\rm GL}_2(F)$ such that $\pi_2 \cong~ ^g\pi_1$. 
Thus, ${\rm GL}_2(F)$ acts transitively on the set of isomorphism
classes of automorphic 
representations
of ${\rm SL}_2$ in a given $L$-packet.
\end{lemma}

\begin{proof}
Since $\pi_1$ is  cuspidal automorphic, it has a  
Whittaker model 
with respect to a character $\psi_1:{\Bbb A}_F/F \longrightarrow {\Bbb C}^*$;
similarly assume that $\pi_2$  has a  
Whittaker model 
with respect to a character $\psi_2:{\Bbb A}_F/F \longrightarrow {\Bbb C}^*$.
As is well known, any two non-trivial characters of ${\Bbb A}_F/F$ differ by 
a scaling from $F^*$, i.e., there exists 
$f \in F^*$  such that 
$\psi_2(x)=\psi_1(fx)$ for all $x \in {\Bbb A}_F/F$.
From the uniqueness of the Whittaker model with respect to ${\rm GL}_2(F_v)$, 
it follows
 that if $\pi_{1,v}$ has a Whittaker model with respect to $\psi_{1,v}$ and 
$\pi_{2,v}$ for $\psi_{2,v}$, and if $\psi_{2,v}(x)=\psi_{1,v}(f_vx)$ for some 
$f_v \in F_v$, then $\pi_{2,v} \cong~^{g}\pi_{1,v}$ where $g$ is any element
of ${\rm GL}_2(F_v)$ with $\det g = f_v$. This completes
the proof of the lemma.
\end{proof}

\begin{cor}\label{cor2.2}
Let $\pi$ be an irreducible representation of ${\rm SL}_2({\Bbb A}_F)$ contained
in  the restriction of a cuspidal
 automorphic representation $\tilde{\pi}$ of ${\rm GL}_2({\Bbb A}_F)$. Then 
$\pi$ is automorphic if and only if $\pi$ has a Whittaker model with
respect to a non-trivial character $\psi: {\Bbb A}_F/F \longrightarrow 
{\Bbb C}^*$.
\end{cor}

\begin{proof} 
Clearly if $\pi$ is cuspidal automorphic, it has a Whittaker model.
Conversely, fix $\pi_1$ to be an automorphic representation of ${\rm SL}_2({\Bbb A}_F)$
contained in $\tilde{\pi}$, and suppose that $\pi_1$ has a Whittaker model 
with respect to a non-trivial character $\psi_1: {\Bbb A}_F/F \longrightarrow {\Bbb C}^*$.
Since the set of non-trivial characters of ${\Bbb A}_F/F$ is parametrized
by $F^*$ as before, there exists $a \in F^*$ such that 
$\psi(x)=\psi_1(ax)$ for all $x \in {\Bbb A}_F/F$. Since ${\rm GL}_2(F)$ operates 
on the set of automorphic representations of ${\rm SL}_2({\Bbb A}_F)$ contained in 
$\tilde{\pi}$, by conjugating $\pi_1$ by 
$\left(
\begin{array}{cc}
a & 0\\
0 & 1
\end{array}
\right)$,
we can assume that $\psi=\psi_1$. By the uniqueness of 
Whittaker model (for ${\rm GL}_2$) $\pi \cong \pi_1$, hence $\pi$ is automorphic.
\end{proof}

\begin{cor}\label{cor2.3}
For an irreducible representation $\pi$ of ${\rm SL}_2({\Bbb A}_F)$ contained in a cuspidal
automorphic representation $\tilde{\pi}$ of ${\rm GL}_2({\Bbb A}_F)$, the following
are equivalent.
\begin{enumerate}
\item $\pi$ has an 
abstract Whittaker model with respect to a character $\psi: {\Bbb A}_F/F 
\longrightarrow {\Bbb C}^*$.
\item  $\pi$ has a nonzero Fourier coefficient with 
respect to $\psi: {\Bbb A}_F/F \longrightarrow {\Bbb C}^*$.
\item $\pi$ is automorphic.
\end{enumerate}
\end{cor}

\begin{lemma}\label{lemma_loc}
Let $k$ be a local field, $\tilde{\pi}$ an irreducible admissible
representation of $GL_n(k)$, and $\pi$ an irreducible subrepresentation 
of $\tilde{\pi}$ restricted to $SL_n(k)$. Let
$$X_{\tilde{\pi}}=\{\chi:k^*\rightarrow {\Bbb C}^* \mid \tilde{\pi}\otimes\chi
\cong \tilde{\pi}\}$$
and
$$G_\pi=\{g\in GL_n(k) \mid \pi^g \cong \pi\}.$$
Then $G_\pi=\displaystyle{\bigcap_{\chi \in X_{\tilde{\pi}}}}{\rm ker}(\chi)$ 
where for a character $\chi$ of $k^*$, 
$${\rm ker}(\chi)=\{g\in GL_n(k) \mid  \chi(\det g) =1\}.$$

\end{lemma}

\begin{proof}
See, for example,
Theorem 4.2 of \cite{gelbart} for a proof of this well-known lemma.
\end{proof}

\begin{lemma}\label{lemma2.4}
If $\pi=\otimes \pi_v$ is an automorphic representation of ${\rm SL}_2({\Bbb A}_F)$
which is not a monomial automorphic representation, then any $\pi^\prime=\otimes \pi^\prime_v$
with $\pi^\prime_v$ in the $L$-packet containing
 $\pi_v$ and equal to $\pi_v$ at almost all places $v$ of $F$ is automorphic.
\end{lemma}

\begin{proof} 
It suffices to prove that $\pi^\prime \cong ~^g\pi$ for
$g\in {\rm GL}_2(F)$. Define 
$$G_\pi=\{g \in {\rm GL}_2({\Bbb A}_F) \mid {^g\pi}\cong \pi\}.$$
Clearly $G_\pi$ contains ${\rm SL}_2({\Bbb A}_F)$ as well as
 ${\Bbb A}^*_F$ embedded
diagonally in ${\rm GL}_2({\Bbb A}_F)$. To prove the lemma, it 
suffices to prove that the double coset 
$${\rm GL}_2(F)\backslash {\rm GL}_2({\Bbb A}_F)/G_\pi,$$
consists of a single element. This is clearly an abelian group 
which is a quotient of ${\Bbb A}^*_F/F^*$. We will prove that this group
is trivial by proving that it has no nontrivial characters.

Assume that $\tilde{\pi} = \otimes \tilde{\pi}_v$ is an irreducible
automorphic representation of ${\rm GL}_2({\Bbb A}_F)$ 
containing the automorphic
representation $\pi$ of ${\rm SL}_2({\Bbb A}_F)$. From Lemma \ref{lemma_loc},
we find that
$$G_\pi=
\left \{ g \in {\rm GL}_2({\Bbb A}_F) \mid \chi(\det g) = 1 ~\forall 
\chi:{\Bbb A}^*_F \rightarrow {\Bbb C}^* {\rm such~ that} ~ \pi \otimes \chi \cong \pi 
\right \}.$$
Therefore the characters of ${\rm GL}_2(F)\backslash {\rm GL}_2({\Bbb A}_F)/G_\pi$
are Gr\"ossencharacters $\chi$ such that $\tilde{\pi} \otimes \chi \cong 
\tilde{\pi}$.
However as $\tilde{\pi}$ is non-monomial, there are no such characters, 
proving the lemma.
\end{proof}

\begin{remark}
The same proof yields that in the monomial case there are either two 
or four orbits of the ${\rm GL}_2(F)$ action on representations of 
${\rm SL}_2({\Bbb A}_F)$ (not necessarily automorphic) belonging to one $L$-packet. 
Exactly one orbit consists of automorphic representations, 
and the other (one or three) orbits do
not have any automorphic representation.
\end{remark}

\begin{remark}\label{remark2.7}
Our proof works more generally for $SL_n({\Bbb A}_F)$ to prove stability of 
primitive representations of $GL_n({\Bbb A}_F) $, i.e., those
automorphic representations $\tilde{\pi} $ 
of $GL_n({\Bbb A}_F) $ for which there are no nontrivial
characters $\chi: {\Bbb A}_F^*/F^* \rightarrow {\Bbb C}^*$ 
with $\tilde{\pi}\otimes\chi \cong \tilde{\pi}$.
\end{remark}

\subsection{Size of $L$-packets}

Lemma \ref{lemma2.4} says that for a global automorphic $L$-packet on
${\rm SL}_2$ which by definition is made up of local packets, 
one can change any local component in its $L$-packet in the
non-monomial case. This
brings us to the interesting question whether the size of a 
non-monomial global
$L$-packets is finite or infinite. This does not seem to have
been studied in the literature, either for ${\rm SL}_2$, or for
other groups. We take this opportunity to make a remark about it.

Observe that since most local components of an automorphic form on
${\rm GL}_2$ are unramified principal series, therefore given by a pair of
complex numbers $(\alpha_v, \beta_v)$, $v$ running over all
but finitely many places of $F$, the question amounts to
whether for infinitely many places $v$ of $F$, the corresponding
principal series representation of
${\rm GL}_2(F_v)$ reduces into more than one component when restricted to 
${\rm SL}_2(F_v)$. This is the case if and only if $\alpha_v = -\beta_v$,
i.e., $\alpha_v+ \beta_v =0$. Thus for modular forms for 
$\Gamma_1(N) \subset {\rm SL}_2({\Bbb Z})$, given by classical Fourier expansion
$$f(z) = \sum_n a_n e^{2\pi i nz},$$
the question amounts to whether $a_p=0$, for infinitely many primes $p$.
As is well known, N. Elkies proved the existence of infinitely many 
such primes,
called {\it supersingular} primes, when the modular form comes from an
elliptic curve. However, existence of infinitely many such primes
is perhaps a special feature of modular forms of weight 2 whose
Fourier coefficients lie in ${\Bbb Z}$; it is not clear what 
to expect for forms of 
higher weight, or for forms whose Fourier coefficients do not lie in 
${\Bbb Z}$ (always of course 
in the non-monomial case). We refer to the article of Kumar Murty, which 
establishes upper bounds for such primes in \cite{murty}. 

\section{Global distinguishedness of an $L$-packet for ${\rm SL}_2$}

We introduce some notation. For any number field $F$, let
$${\Bbb A}^1_F=\{x=(x_v)\in {\Bbb A}^*_F \mid |x|=\prod_v |x_v|_v=1\}.$$
By the product formula, $F^* \subseteq {\Bbb A}^1_F$, and it is well-known that
$F^*\backslash {\Bbb A}^1_F$ is a compact group.

\indent

Similarly, let
$${\rm GL}_2^1({\Bbb A}_F)=\{g\in {\rm GL}_2({\Bbb A}_F)\mid \det g \in {\Bbb A}^1_F\}.$$

\begin{lemma}\label{lemma3.1}
Let $E$ be a quadratic extension of a number field $F$. Let $\phi$ be a cusp 
form on ${\rm GL}_2({\Bbb A}_E)$ whose central character restricted to ${\Bbb A}^*_F$
is trivial. Then
$${\rm vol}(F^*\backslash {\Bbb A}^1_F) \cdot 
\int_{{\Bbb A}^*_F {\rm GL}_2(F)\backslash {\rm GL}_2({\Bbb A}_F)}\phi(g)dg=
\int_{{\rm GL}_2(F)\backslash {\rm GL}_2^1({\Bbb A}_F)}
\phi(g)dg.$$
\end{lemma}

\begin{proof}
 The absolute convergence of the two integrals above is a  
consequence of the fact that cusp forms are bounded, and that we are 
dealing with spaces with finite volume.
The equality of the integrals is clear as the natural group homomorphism 
from ${\rm GL}_2^1({\Bbb A}_F)$ to ${\rm PGL}_2({\Bbb A}_F)$ is surjective with kernel
consisting of $x \in {\Bbb A}^*_F$ with $| x|^2=1$ which is nothing but 
${\Bbb A}^1_F$.
\end{proof}

\begin{prop}\label{prop3.2}
Let $E$ be a quadratic extension of a number field $F$. Let $\phi$ be a cusp 
form on ${\rm GL}_2({\Bbb A}_E)$. Then
$$\int_{{\rm SL}_2(F)\backslash {\rm SL}_2({\Bbb A}_F)}\phi(g)dg = 
\frac{1}{{\rm vol}(F^*\backslash {\Bbb A}^1_F)}\sum_\eta \int_{
{\rm GL}_2(F)\backslash {\rm GL}_2^1({\Bbb A}_F)}\phi(g)\eta(\det g)dg$$
where the sum on the right hand side of the equality sign is over all 
characters $\eta$ of the compact abelian group $F^*\backslash {\Bbb A}^1_F$.
\end{prop}

\begin{proof} 
We note that for a locally compact topological group 
$G$ with closed subgroups $H_1 \subset H_2$, which are all 
assumed to be unimodular, there exists
a choice of invariant measures on $H_1 \backslash G$, $H_2 \backslash G$, 
$H_1 \backslash H_2$, denoted by $d_1g,d_2g,dh$,
such that for a function $f \in L^1(H_1\backslash G),$

$$ \int_{H_1\backslash G} f(g)d_1g = \int_{ H_2 \backslash G}\left ( 
\int_{H_1 \backslash H_2} f(hg)dh \right )d_2g.$$

Applying this general result to  ${\rm GL}_2(F) \subset {\rm GL}_2(F) {\rm SL}_2({\Bbb A}_F) 
\subset {\rm GL}_2^1({\Bbb A}_F)$, we have,

\begin{align}
\int_{{\rm GL}_2(F)\backslash {\rm GL}_2^1({\Bbb A}_F)}\phi(g)dg=
\int_{F^*\backslash{\Bbb A}^1_F}\left(\int_{{\rm SL}_2(F)\backslash {\rm SL}_2(
{\Bbb A}_F)}\phi(gx)dg \right) dx.
\end{align}

\noindent
Define
$$F(x)=\int_{{\rm SL}_2(F)\backslash {\rm SL}_2({\Bbb A}_F)}\phi(gx)dg$$
for $x\in {\Bbb A}^1_F$ embedded inside ${\rm GL}_2^1({\Bbb A}_F)$ as
$\left(\begin{array}{cc} x & 0 \\ 0 & 1 \end{array} \right)$.
Clearly $F(x)$ is a function on $F^*\backslash {\Bbb A}^1_F$.
By the Fourier inversion theorem
$$F(1)=\frac{1}{{\rm vol}(F^*\backslash {\Bbb A}^1_F)}\sum_\eta \int_{F^*\backslash {\Bbb A}^1_F}F(x)\eta(x)dx,$$
where the sum on the right hand side of the equality sign is over all 
characters $\eta$ of the compact abelian group $F^*\backslash {\Bbb A}^1_F$.
Thus by (1), the proof of the proposition is completed.
\end{proof}

We next note the following lemma.

\begin{lemma}
For a character $\chi: {\Bbb A}^*_F/F^* \longrightarrow {\Bbb C}^*$, there
exists a character $ \eta: {\Bbb A}^*_F/F^* \longrightarrow {\Bbb C}^*$ 
such that $\chi = \eta^2$ if and only if there are no local obstruction
to solving $\chi  = \eta^2$, i.e., if $\chi = \prod \chi_v$, then
$\chi_v(-1) = 1$ for all places $v$ of $F$.
\end{lemma}
\begin{proof} The proof follows easily by analyzing the
exact sequence of topological abelian groups 
$$0 \rightarrow A[2] \rightarrow A \rightarrow A,$$
with $A =  {\Bbb A}^*_F/F^*,$ and $A[2] = \{a \in A| a^2 =1 \}$ 
together with the fact that an element of $F^*$ 
is a square if and only if it is a square in $F_v^*$ for all
places $v$ of $F$.
\end{proof}

\begin{prop}\label{prop3.3}
If $\tilde{\pi}$ is a cusp form on ${\rm GL}_2({\Bbb A}_E)$ which is 
distinguished by 
${\rm SL}_2({\Bbb A}_F)$, then there is a Gr\"ossencharacter $\eta$ of $F^*\backslash
{\Bbb A}^*_F$ such that $\tilde{\pi}$ is $\eta$-distinguished 
for ${\rm GL}_2({\Bbb A}_F)$.
Conversely if $\tilde{\pi}$ is $\eta$-distinguished for 
some Gr\"ossencharacter $\eta$
of $F^*\backslash {\Bbb A}^*_F$, then $\tilde{\pi}$ is 
${\rm SL}_2({\Bbb A}_F)$- 
distinguished. Hence there is a member of the $L$-packet of automorphic
representations of 
${\rm SL}_2({\Bbb A}_E)$ determined by $\tilde{\pi}$ which is globally
${\rm SL}_2({\Bbb A}_F)$-distinguished.
\end{prop} 

\begin{proof} 
As $\tilde{\pi}$ is distinguished by 
${\rm SL}_2({\Bbb A}_F)$, it is locally distinguished. Hence  the
central character $\omega_{\tilde{\pi}}$ of $\tilde{\pi}$ takes the value 1
at $-1$
locally at all places $v$ of $F$. Therefore by the
previous lemma, we can  
assume that $\omega_{\tilde{\pi}}$ restricted to ${\Bbb A}^*_F$ 
is the square  of a 
Gr\"ossencharacter on ${\Bbb A}^*_F$ and hence by twisting that   
the central character of $\tilde{\pi}$ restricted to ${\Bbb A}^*_F$
is trivial. (Actually, by the same 
argument $\omega_{\tilde{\pi}}$ itself 
is the square  of a 
Gr\"ossencharacter on ${\Bbb A}^*_E$ and hence by twisting we can assume 
that    the central character of $\tilde{\pi}$ 
is trivial, but this is not relevant for us.) 

Now combining Lemma \ref{lemma3.1} and Proposition \ref{prop3.2}, 
and assuming that ${\rm vol}(F^*\backslash {\Bbb A}^1_F)= 1$,
we have:
\begin{eqnarray*}
\int_{{\rm SL}_2(F)\backslash {\rm SL}_2({\Bbb A}_F)}\phi(g)dg 
& = & {\sum_{\eta: F^*\backslash 
{\Bbb A}^1_F \rightarrow {\Bbb C}^*}} 
\int_{{\rm GL}_2(F)\backslash {\rm GL}_2^1({\Bbb A}_F)}\phi(g)
\eta(\det g)dg \\
& = &  {\sum_{ 
\begin{array}{c} \eta: F^*\backslash 
{\Bbb A}^1_F \rightarrow {\Bbb C}^* \\
\eta^2 =1 \end{array}}} 
\int_{{\rm GL}_2(F)\backslash {\rm GL}_2^1({\Bbb A}_F)}\phi(g)
\eta(\det g)dg \\
& = &   {\sum_{
\begin{array}{c} \tilde{\eta}: F^*\backslash 
{\Bbb A}^*_F \rightarrow {\Bbb C}^* \\
\tilde{\eta}^2 =1 \end{array}}}
\int_{{\Bbb A}^*_F {\rm GL}_2(F)
\backslash {\rm GL}_2({\Bbb A}_F)}
\phi(g)\tilde{\eta}(\det g)dg
\end{eqnarray*}

Thus if $\tilde{\pi}$ is distinguished by ${\rm SL}_2({\Bbb A}_F)$, 
then it is $\eta$ 
distinguished by ${\rm GL}_2({\Bbb A}_F)$ for some Gr\"ossencharacter $\eta$ of
$F^*\backslash {\Bbb A}^*_F$. 

Conversely, assume that $\tilde{\pi}$ is $\eta$ distinguished by 
${\rm GL}_2({\Bbb A}_F)$, 
and \\
$\int_{{\rm SL}_2(F)\backslash {\rm SL}_2({\Bbb A}_F)}
\phi(g)dg =0$ for all $\phi \in \tilde{\pi}$. Twisting by a character, we 
assume that $\eta=1$. Then, in particular,
 $\int_{{\rm SL}_2(F)\backslash {\rm SL}_2({\Bbb A}_F)}
\phi(gx)dg =0$ for all $x \in {\rm GL}_2( {\Bbb A}_F)$.
By the identity $(1)$ in the proof of Proposition \ref{prop3.2}, we get,
$$\int_{{\rm GL}_2(F)\backslash {\rm GL}_2^1({\Bbb A}_F)}\phi(g)dg=
\int_{F^*\backslash{\Bbb A}^1_F}
\left(
\int_{{\rm SL}_2(F)\backslash {\rm SL}_2(
{\Bbb A}_F)}\phi(gx)dg 
\right) dx=0$$
which, by an application of Lemma \ref{lemma3.1}, 
is a contradiction to $\tilde{\pi}$ being 
distinguished by ${\rm GL}_2({\Bbb A}_F)$, completing
the proof of Proposition \ref{prop3.3}. 
\end{proof}

\begin{remark} 
When we talk of $\chi$ distinguished representation,
$\chi$ is a character of $F^*$ or ${\Bbb A}^*_F/F^*$ as the case may be,
whereas in many calculations, we have to extend this character to 
a character of 
$E^*$ or $ {\Bbb A}^*_E/E^*$ which we often continue to write as $\chi$.
The end results naturally depend only on $\chi$ on $F^*$ or 
${\Bbb A}^*_F/F^*$, and not on the extension chosen. 
\end{remark}
                
\section{Criterion for global distinguishedness for ${\rm SL}_2$}

We begin with the following local result which follows from Theorem 1.1 of
\cite{anand}.

\begin{lemma}\label{lemma4.1}
Let $K$ be a quadratic extension of a local field $k$. Let $\pi$ be an 
irreducible admissible representation of ${\rm SL}_2(K)$ 
contained in an irreducible
admissible representation $\tilde{\pi}$ of ${\rm GL}_2(K)$ which is
distinguished by $GL_2(k)$. Then $\pi$ is
distinguished by ${\rm SL}_2(k)$ if and only if $\pi$ has a Whittaker model
with respect to a character of $K$ which is trivial on $k$.
\end{lemma}

Here is the theorem about global distiguishedness of an
automorphic representation of ${\rm SL}_2({\Bbb A}_E)$ which is the global
analogue of the local result contained in Lemma \ref{lemma4.1}.

\begin{thm}\label{thm4.2}
Let $\pi$ be an automorphic representation of
${\rm SL}_2({\Bbb A}_E)$ contained in a cuspidal
 automorphic representation
$\tilde{\pi}$ of ${\rm GL}_2({\Bbb A}_E)$. Suppose that $\tilde{\pi}$
is distinguished by ${\rm GL}_2({\Bbb A}_F)$. Then $\pi$ is distinguished
by ${\rm SL}_2({\Bbb A}_F)$ if and only if it has a Whittaker model
with respect to 
 a non-trivial character of ${\Bbb A}_E/E$ trivial on ${\Bbb A}_F/F$.
\end{thm}

The proof of this theorem will use the following lemma.
\begin{lemma}\label{lemma4.3}
Let $\phi$ be a square integrable function on ${\rm SL}_2(F)\backslash {\rm SL}_2({\Bbb A}_F)$ such that
$$\int_{N(F)\backslash N({\Bbb A}_F)}\phi(ng)dn=0$$
for all $g\in {\rm SL}_2({\Bbb A}_F)$ where $N$ is the group of all upper triangular
unipotent matrices in ${\rm SL}_2$. Then $$\int_{{\rm SL}_2(F)\backslash 
{\rm SL}_2({\Bbb A}_F)}\phi(g)dg=0.$$
\end{lemma}

\begin{proof} 
The condition on $\phi$ implies that it is a cusp form, hence 
it belongs to the (completion) of the direct sum of cuspidal automorphic 
representations in 
$L^2({\rm SL}_2(F)\backslash {\rm SL}_2({\Bbb A}_F))$. The integral
$$f \mapsto
\int_{{\rm SL}_2(F)\backslash {\rm SL}_2({\Bbb A}_F)}f(g)dg$$
is an ${\rm SL}_2({\Bbb A}_F)$-invariant linear form, 
and hence must be trivial 
on any irreducible representation which is not trivial, hence on any 
irreducible cuspidal representation, and therefore on their sum too. It follows
that $\int_{{\rm SL}_2(F)\backslash {\rm SL}_2({\Bbb A}_F)}\phi(g)dg=0$.
\end{proof}

\begin{proof}[Proof of Theorem~\ref{thm4.2}]
Suppose that $\pi$ is distinguished by 
${\rm SL}_2({\Bbb A}_F)$. Then 
$\int_{{\rm SL}_2(F)\backslash {\rm SL}_2({\Bbb A}_F)}\phi(g)dg 
\neq 0$, for some $\phi \in \pi$. By the previous lemma, this implies that
$\int_{N(F)\backslash N({\Bbb A}_F)}\phi(ng)dn\neq 0$ for some $g \in 
{\rm SL}_2({\Bbb A}_F)$. By considering the translate of $\phi$ by $g$, 
one can 
in fact assume that $\int_{N(F)\backslash N({\Bbb A}_F)}\phi(n)dn \neq 0$.
Now $\phi$ is a cusp form on ${\rm SL}_2(E)\backslash {\rm SL}_2({\Bbb A}_E)$. 
Considering
it as a function on $N(E)\backslash N({\Bbb A}_E)$, which we henceforth write
as ${\Bbb A}_E/E$, and expanding it as a Fourier series, we have 
$$\phi(n)=\displaystyle{\sum_\psi}\hat{\phi}(\psi)\psi(n)$$
where $\psi$ runs over all characters $\psi:{\Bbb A}_E/E \longrightarrow 
{\Bbb C}^*$, and $$\hat{\phi}(\psi)=\int_{{\Bbb A}_E/E}\phi(v)\psi(-v)dv.$$
Since the integral of a non-trivial character on ${\Bbb A}_F/F$ is zero, we
find that
$$\int_{N(F)\backslash N({\Bbb A}_F)}\phi(n)dn=\sum_{\psi}\hat{\phi}(\psi),$$
where $\psi$ runs over all characters $\psi:{\Bbb A}_E/E\longrightarrow
{\Bbb C}^*$ which are trivial on ${\Bbb A}_F$. Since $\int_{N(F)\backslash 
N({\Bbb A}_F)}\phi(n)dn \neq 0$, 
there must be a $\psi$ which is trivial on ${\Bbb A}_F/F$, for which 
$\hat{\phi}(\psi)\neq 0$. By the cuspidality condition, $\psi$ must be 
non-trivial. This proves the existence of a Whittaker model with
respect to a character of ${\Bbb A}_E/E$ trivial on ${\Bbb A}_F/F$.

We now prove the converse statement, i.e., if $\pi$ has a Whittaker model with respect to
a character $\psi:{\Bbb A}_E/E \longrightarrow {\Bbb C}^*$ which is trivial on
${\Bbb A}_F/F$, then $\pi$ is distinguished. For this observe that by 
Proposition \ref{prop3.3}, $\tilde{\pi}$ is ${\rm SL}_2({\Bbb A}_F)$-distinguished, 
and hence
some cuspidal representation in the global $L$-packet of $\pi$ is
${\rm SL}_2({\Bbb A}_F)$-distinguished. By Lemma \ref{lemma2.1}, 
we can assume that $^g\pi$ is
distinguished for some $g\in {\rm GL}_2(E)$, hence from what has been just proved,
$^g\pi$ has a Whittaker model by a character $\psi^\prime:{\Bbb A}_E/E
\longrightarrow {\Bbb C}^*$ which is trivial on ${\Bbb A}_F/F$. But we are
given that $\pi$ has a Whittaker model by a character $\psi:{\Bbb A}_E/E
\longrightarrow {\Bbb C}^*$ which is trivial on ${\Bbb A}_F/F$. Since the set
of non-trivial characters of ${\Bbb A}_E/E$ trivial on 
${\Bbb A}_F/F$ is a principal homogeneous space for $F^*$, and since clearly 
$\pi$ is distinguished by ${\rm SL}_2({\Bbb A}_F)$ if and only is $^h\pi$ is for
any $h\in {\rm GL}_2(F)$, we can assume that $\psi=\psi^\prime$, i.e., both $\pi$ and
$^g\pi$ have Whittaker models by the same character $\psi$ and 
$^g\pi$ is distinguished.
But by the uniqueness of Whittaker model (for ${\rm GL}_2$), this implies that
$\pi={^g\pi}$, and hence $\pi$ is distinguished.
\end{proof}

\noindent
\begin{remark} 
E. Lapid has pointed out to us that Lemma \ref{lemma4.3}
can also be proved as 
follows. Every cusp form is orthogonal to any pseudo-Eisenstein
series, and the pseudo-Eisenstein
series contain the constant functions in their closure, thus a cusp form is 
orthogonal to the constants.
\end{remark}

\section{Locally but not globally distinguished I}

In this section we use the theorem of the previous 
section to  show that 
there are cuspidal representations
of ${\rm SL}_2({\Bbb A}_E)$ which are not distinguished by ${\rm SL}_2({\Bbb A}_F)$
but for which each of its 
local component is ${\rm SL}_2$-distinguished. To this end, 
fix a nontrivial character $\psi: {\Bbb A}_E/E \longrightarrow {\Bbb C}^*$
which is trivial on ${\Bbb A}_F$. Let $\pi$
be a cuspidal representation of ${\rm SL}_2({\Bbb A}_E)$ which occurs in the 
restriction of a cuspidal representation $\tilde{\pi}$ of ${\rm GL}_2({\Bbb A}_E)$
which is distinguished by ${\rm GL}_2({\Bbb A}_F)$.

Our examples will depend on understanding and identifying the distinguished
parts of the restriction of $\tilde{\pi}$ to the successive subgroups
$$GL_2({\Bbb A}_E) \supseteq  
{\Bbb A}^*_E {\rm SL}_2({\Bbb A}_E){\rm GL}_2({\Bbb A}_F) \supseteq
{\Bbb A}^*_E {\rm SL}_2({\Bbb A}_E) {\rm GL}_2(F) \supseteq
SL_2({\Bbb A}_E).$$

We denote by $\pi^\prime$ the irreducible 
representation of $G^\prime={\Bbb A}^*_E {\rm SL}_2({\Bbb A}_E) {\rm GL}_2({\Bbb A}_F)$ 
that occurs in the restriction of $\tilde{\pi}$ to $G^\prime$,
and  which is 
$\psi$-generic. 
By Lemma \ref{lemma4.1}, $\pi^\prime$ is the unique
irreducible component of the restriction of $\tilde{\pi}$ to $G^\prime$
which is abstractly distinguished by ${\rm SL}_2({\Bbb A}_F)$.  
Further,  an irreducible representation of 
${\rm SL}_2({\Bbb A}_E)$ occurring in $\tilde{\pi}$ is abstractly 
distinguished with respect to 
${\rm SL}_2({\Bbb A}_F)$ if and only if it occurs in the restriction of $\pi^\prime$ 
to ${\rm SL}_2({\Bbb A}_E)$. From Theorem \ref{thm4.2}, it follows that there is exactly
one irreducible 
cuspidal representation of $G^{\prime\prime}= 
{\Bbb A}^*_E {\rm SL}_2({\Bbb A}_E) {\rm GL}_2(F)$
occurring in the space of functions in $\tilde{\pi}$ 
that is distinguished by ${\rm SL}_2({\Bbb A}_F)$, say $\pi^{\prime\prime}$.
Also, an irreducible cuspidal representation of 
${\rm SL}_2({\Bbb A}_E)$ occurring in the space of functions in 
$\tilde{\pi}$ is distinguished by 
${\rm SL}_2({\Bbb A}_F)$ if and only if it occurs in the restriction of 
$\pi^{\prime\prime}$ to ${\rm SL}_2({\Bbb A}_E)$. Now if we choose $\tilde{\pi}$
such that $\tilde{\pi} \otimes \omega \cong \tilde{\pi}$ where $\omega$
is a character of ${\Bbb A}_E^*/E^*$ with 
 non-trivial restriction to ${\Bbb A}^*_F$, 
then $\pi^\prime \otimes \omega \cong \pi^\prime$, since $\pi^\prime$
is the unique irreducible representation of $G^\prime$ which is 
(abstractly) $\psi$-generic.
Hence the restriction of
$\pi^\prime$ to $G^{\prime\prime}$ is not irreducible. Hence we get 
cuspidal 
representations of ${\rm SL}_2({\Bbb A}_E)$ which appear in the restriction of 
$\pi^\prime$ but not in the restriction of $\pi^{\prime\prime}$. These 
representations are abstractly distinguished but not distinguished.

It remains to construct cuspidal representations $\tilde{\pi}$ of 
${\rm GL}_2({\Bbb A}_E)$ which are distinguished by ${\rm GL}_2({\Bbb A}_F)$
such that $\tilde{\pi} \otimes \omega \cong \tilde{\pi}$ where
$\omega$ restricts non-trivially to ${\Bbb A}^*_F$. 
We need the following lemma.
\begin{lemma}\label{lemma4.4}
Let $L/F$ be a quadratic extension of number fields. Given a positive
integer $n$, there exists a Gr\"ossencharacter $\eta$ of ${\Bbb A}^*_L$
of order $n$ such that $\eta$ has trivial restriction to ${\Bbb A}^*_F$.
\end{lemma}

\begin{proof} 
Let $v$ be a place of $F$ that splits in $L$, say $v=w_1w_2$,
such that $F_v$ has odd residue characteristic. Let $\eta_{w_1}$ be
a character of order $n$ of $L^*_{w_1}$. Consider the character
$(\eta_{w_1},1)$ of $L^*_{w_1} \times L^*_{w_2}$. By Grunwald-Wang theorem,
we get a Gr\"ossencharacter of ${\Bbb A}^*_L$ of order $n$ whose component
at the place above $v$ is $(\eta_{w_1},1)$. It follows from our construction 
that $\eta/\eta^\tau$
is also a Gr\"ossencharacter of order $n$, where $\tau$ is the non-trivial
element of $Gal(L/F)$. Further $\eta/\eta^\tau$ restricted to ${\Bbb A}^*_F$
is trivial.
\end{proof}

Now let $\eta$ be a Gr\"ossencharacter of ${\Bbb A}^*_L$ of order 8
such that $\eta$ has trivial restriction to ${\Bbb A}^*_F$.
Let $M$ be the quadratic extension of $L$ such that $\eta^4=\omega_{M/L}$.
Since $\omega_{M/L}$ has trivial restriction to ${\Bbb A}^*_F$, we see that
there is a quadratic extension $E$ of $F$ such that $M=EL$ (cf. Corollary
\ref{cor6.2}).
The conditions on $\eta$ imply that the representations 
${\rm Ind}_{W_L}^{W_F}\eta$ and ${\rm Ind}_{W_L}
^{W_F}\eta^2$ of $W_F$, the Weil group of $F$, 
are irreducible and the restriction of   
${\rm Ind}_{W_L}^{W_F}\eta^2$ to $W_E$, the Weil group of $E$,
is a sum of two distinct Gr\"ossencharacters, necessarily of the form
$\gamma$ and $\gamma^\sigma$. Now let $\rho_{\tilde{\pi}}$
be the restriction to $W_E$ of the representation 
${\rm Ind}_{W_L}^{W_F}\eta$ of $W_F$,
and $\tilde{\pi}$ the associated automorphic form
on ${\rm GL}_2({\Bbb A}_E)$. Since $\eta^4 \neq 1$,
$\rho_{\tilde{\pi}}$ is an irreducible representation. 

Let $r(\rho_{\tilde{\pi}})$ be the 4 dimensional 
representation of $W_F$ obtained from  $\rho_{\tilde{\pi}}$ 
of $W_E$ 
by the process of twisted tensor induction. 
It is a general and simple fact that if $H$ is a subgroup of a group $G$
of index two and $V$ a representation of $G$, then
$$r(V|_H) \cong {\rm Sym}^2(V)
\oplus \wedge^2(V)\cdot 
\omega_{G/H}$$
where $\omega_{G/H}$ is the nontrivial character of $G$ trivial on $H$.
Applying this to our situation, we have:

\begin{align*}
r(\rho_{\tilde{\pi}}) &\cong {\rm Sym}^2({\rm Ind}_{W_L}^{W_F}\eta)
\oplus \wedge^2({\rm Ind}_{W_L}^{W_F}\eta)\cdot 
\omega_{E/F}\\
 &\cong {\rm Ind}_{W_L}^{W_F}\eta^2 \oplus 1 \oplus
\omega_{E/F}\omega_{L/F}.
\end{align*} 

Since $r(\rho_{\tilde{\pi}})$ contains the trivial representation, it 
follows from the known theorems, as recalled in the introduction, 
that $\tilde{\pi}$ is distinguished by ${\rm GL}_2({\Bbb A}_F)$.
Also:

\begin{align*}
\rho_{\tilde{\pi}} \otimes \rho_{\tilde{\pi}^\sigma}
 \cong \rho_{\tilde{\pi}} \otimes 
\rho_{\tilde{\pi}}^\vee = 1 \oplus \omega_{L/F} \circ Nm_{E/F} \oplus
\gamma \oplus \gamma^\sigma.
\end{align*}

Therefore $\gamma$ is a self-twist for $\tilde{\pi}$. Observe that
$\gamma$ has non-trivial restriction to ${\Bbb A}^*_F$. We have thus
proved the following theorem.

\begin{theorem} 
There is a  cuspidal representation of ${\rm SL}_2({\Bbb A}_E)$ which is 
 not distinguished by ${\rm SL}_2({\Bbb A}_F)$
but for which each of its 
local component is ${\rm SL}_2(F_v)$-distinguished. 
\end{theorem}

The above analysis also gives the following proposition.
\begin{prop}\label{prop4.5}
Let $\tilde{\pi}$ be a non-monomial cuspidal representation of ${\rm GL}_2({\Bbb A}_E)$ 
that is distinguished by ${\rm GL}_2({\Bbb A}_F)$.
Then any irreducible cuspidal representation 
of ${\rm SL}_2({\Bbb A}_E)$  in the $L$-packet associated to 
$\tilde{\pi}$ that is abstractly distinguished with respect to 
${\rm SL}_2({\Bbb A}_F)$ is in fact distinguished by ${\rm SL}_2({\Bbb A}_F)$.
\end{prop}

\begin{proof}
Note that since $\tilde{\pi}$ is non-monomial, it cannot be 
$\chi$-distinguished with respect to ${\rm GL}_2({\Bbb A}_F)$ for any non-trivial
Gr\"ossencharacter $\chi$ of ${\Bbb A}^*_F$ (see for example Corollary 
\ref{cor} below). Suppose that $\mu$ is a character
of ${\Bbb A}^*_E$  (not necessarily a 
Gr\"ossencharacter) such that $\tilde{\pi} \otimes \mu \cong \tilde{\pi}$
and such that $\mu$ restricted to $F^*$ is trivial. Since
$r(\tilde{\pi} \otimes \mu)=r(\tilde{\pi})\otimes \mu|_{_{{\Bbb A}^*_F}}$,
it follows that $\tilde{\pi}$ is distinguished with respect to
the Gr\"ossencharacter $\mu|_{_{{\Bbb A}^*_F}}$. This forces 
$\mu|_{_{{\Bbb A}^*_F}}=1$. In other words, any irreducible representation
of $G^\prime$ that occurs in the restriction of $\tilde{\pi}$ to $G^\prime$
restricts irreducibly to $G^{\prime\prime}$. This proves the proposition.
\end{proof}

\section{Factorization}

In this section we analyze whether the period integral on ${\rm SL}_2$ 
is factorizable or not; it is also common to use the word ``Eulerian'' for 
``factorizable''. 
We begin by making a precise definition of factorization
of a linear form $\ell$ on $ \otimes_v \pi_v$, a restricted
direct product of vector spaces $\pi_v$ with respect to 
vectors $w^0_v \in \pi_v$ where $v$ runs over any infinite set, 
say $X$, such as the set of  places of a number field.

We say that $\ell$ is factorizable, if 
there are linear forms $\ell_v$ for each $v \in X$ such that
$\ell_v(w^0_v)=1$ outside a finite subset $T$ of $X$, and such that
for any finite subset $S$ of $X$ containing $T$, 
$$\ell(w_S \otimes w^S) = (\otimes_{v \in S} \ell_v)(w_S),$$
where $  w_S \otimes w^S$ is a vector in $\otimes_{v \in X} \pi_v$
with $w_S \in \otimes_{v\in S} \pi_v$, and $w^S = 
\otimes_{v \not \in S} w^0_v$.

We state the following three elementary lemmas without proof.
\begin{lemma}\label{lemma5.0}
Let $G$ be an algebraic group defined over a number field $F$. Let 
$\pi = \otimes \pi_v$ be an irreducible admissible 
representation of $G({\Bbb A}_F)$. 
Suppose that $H$ is an algebraic subgroup of 
$G$ defined over $F$ such that for each place $v$ of $F$, the irreducible
representation $\pi_v$ of $G(F_v)$ has at most one dimensional
space of $H(F_v)$-invariant forms. Then an $H({\Bbb A}_F)$-invariant
linear form on $\pi$ is factorizable.  
\end{lemma}

\begin{lemma}\label{lemma5.1}
Suppose that $\pi^\prime_v$ is a subspace of $\pi_v$
(containing the vector $w_v$ for almost all $v$), and $\ell$
is a factorizable linear form on $\pi = \otimes_v \pi_v$,
then the restriction of $\ell$ to $\pi^\prime=\otimes_v \pi^\prime_v$
is also factorizable.
\end{lemma}

\begin{lemma}\label{lemma5.2}
Suppose that $\ell_i$ are finitely many
factorizable linear forms $\ell_i = \otimes_v
\ell_{i,v}$ on $\pi= \otimes_v \pi_v$,
then $\ell = \sum_i \ell_i$ is not factorizable
if there is an  infinite subset 
$Y \subset X$ such that the subspace of linear forms on $\pi_v$ 
generated by $
\ell_{i,v}$ has dimension $> 1$ for $v \in Y$.
\end{lemma}

Before we state the main theorem, we prove the following proposition.
This is the global analogue of Proposition 4.2 of \cite{anand}.
\begin{prop}\label{prop}
Let $\pi$ be a cuspidal representation of ${\rm GL}_2({\Bbb A}_E)$
which is (globally) 
distinguished with respect to ${\rm SL}_2({\Bbb A}_F)$. 
Then the sets
$$X=\left\{\chi \in \widehat{{\Bbb A}^*_F/F^*} |
\begin{tabular}{c}
\mbox{ $\pi$ is $\chi$-distinguished}\\ 
\mbox{ with respect to ${\rm GL}_2({\Bbb A}_F)$}
\end{tabular}
\right\}$$ and 
$$Y=\{\mu \in \widehat{{\Bbb A}^*_E/E^*}\mid 
\pi \otimes \mu \cong \pi; \mu|_{_{{\Bbb A}^*_F}}=1\}$$ 
have the same cardinality; in fact 
$\chi \mapsto \chi \circ N_{E/F}$
induces an isomorphism of $X$ onto $Y$ if $\pi$ is 
${\rm GL}_2({\Bbb A}_F)$-distinguished.
\end{prop} 

\begin{proof}
We assume without loss of generality that $\pi$ 
is (globally) distinguished
with respect to ${\rm GL}_2({\Bbb A}_F)$. Then we give explicit 
maps from $X$ to $Y$ and from $Y$ to $X$.

For $\chi \in X$, let $\tilde{\chi}$ be a character of ${\Bbb A}^*_E/E^*$ 
restricting to $\chi$ on ${\Bbb A}^*_F$. 
Then we have $\pi^\vee \cong \pi^\sigma$
and $(\pi \otimes \tilde{\chi}^{-1})^\vee 
\cong (\pi \otimes \tilde{\chi}^{-1})^\sigma$, and therefore
we get
$$\pi \cong \pi \otimes \chi \circ N_{E/F}.$$
Note that since $\pi$ 
is both distinguished and $\chi$-distinguished with respect to 
${\rm GL}_2({\Bbb A}_F)$,
consideration of the central character implies that 
$\omega_\pi|_{_{{\Bbb A}^*_F}}=\omega_\pi|_{_{{\Bbb A}^*_F}}\chi^{-2}=1$. 
Therefore
$\chi^2=1$, thus $\chi \circ N_{E/F} \in Y$.
This allows us to define a
map from $X$ to $Y$ by sending $\chi$ to $\chi \circ N_{E/F}$.

If $\mu \in Y$, then, since $\mu|_{_{{\Bbb A}^*_F}}=1$, and $\mu^2=1$,
we have that $\mu$ factors through the norm map $N_{E/F}$.
Let $\mu=\eta \eta^\sigma$ for a Gr\"ossencharacter $\eta$ of ${\Bbb A}^*_E$.
Now consider the representation $\pi \otimes \eta$.
Observe that 
$(\pi \otimes \eta)^\vee \cong 
(\pi \otimes \eta)^\sigma,$
and that 
$\omega_{_{\pi \otimes \eta}}|_{_{F^*}}=1$.
Therefore
$\pi \otimes \eta$ is either distinguished with respect to ${\rm GL}_2(
{\Bbb A}_F)$ or
$\omega_{_{E/F}}$-distinguished with respect to ${\rm GL}_2({\Bbb A}_F)$.
We map $\mu$ to $\eta|_{_{{\Bbb A}^*_F}}$ or $\eta|_{_{{\Bbb A}^*_F}}
\omega_{_{E/F}}$
accordingly. Clearly the above two maps are inverses of each other
and hence $X$ and $Y$ have the same cardinality, completing the proof
of the proposition.
\end{proof}

\begin{cor}\label{cor}
A non-monomial automorphic representation is $\chi$-distinguished
for at most one Gr\"ossencharacter. A distinguished 
monomial automorphic representation is
$\chi$-distinguished for at least two (and at most four) Gr\"ossencharacters 
$\chi$ of ${\Bbb A}^*_F$.
\end{cor}
\begin{proof}
We need to supply a proof only for monomial representations. Let $\tilde{\pi}$ 
be a distinguished monomial automorphic representation of 
${\rm GL}_2({\Bbb A}_E)$.
We need to show that there exists a non-trivial Gr\"ossencharacter  $\mu$ 
of ${\Bbb A}_E^*/E^*$ with $ \mu|_{_{{\Bbb A}^*_F}}=1$, and
$\tilde{\pi} \otimes \mu \cong \tilde{\pi}$. Since $\tilde{\pi}$ 
is distinguished,
$\tilde{\pi}^\vee \cong \tilde{\pi}^\sigma$ from which it follows that if
$\tilde{\pi} \otimes \mu \cong \tilde{\pi}$, then  
$\tilde{\pi} \otimes \mu^\sigma \cong \tilde{\pi}$ 
also,
and hence  $\tilde{\pi} \otimes (\mu \mu^\sigma) \cong \tilde{\pi}$. 
Since $\tilde{\pi}$ 
is monomial, it
has a non-trivial self-twist $\mu$, hence we are done unless this
self-twist $\mu$ restricted to ${{{\Bbb A}^*_F}}$  equals $\omega_{E/F}$.
But this would mean that $\tilde{\pi}$ is both distinguished, and
 $\omega_{E/F}$-distinguished, which is not possible as observed in the
introduction of the paper (as it would then 
contribute a pole of order 2 to the Rankin product $L$-function 
$L(s, \pi \otimes \pi^\sigma)$).
\end{proof}

From the proof of the previous corollary, we isolate the following fact
which we will have occasion to use in the next theorem about factorization.

\begin{lemma}\label{lemma6.5}
Let $\tilde{\pi}$ 
be a distinguished cuspidal automorphic representation of 
${\rm GL}_2({\Bbb A}_E)$.
Then if $\tilde{\pi} \otimes \mu \cong \tilde{\pi}$,  
$\mu$ restricted to ${{{\Bbb A}^*_F}}$  cannot be equal to $\omega_{E/F}$.
Thus, if $\tilde{\pi}$ is a monomial representation coming from a quadratic extension
$M$ of $E$, $M$ cannot be a cyclic quartic extension of $F$.
\end{lemma}
\begin{proof} The last conclusion is a consequence
 of class field theory, see
Corollary \ref{cor6.2} below.
\end{proof}

Before we proceed further, we note the following lemma from 
class field theory.
\begin{lemma}\label{lemma6.1}
Let $E$ be a finite extension of a number field or a local field $F$.
Let $\chi: {\Bbb A}^*_F/F^* \longrightarrow {\Bbb C}^*$ (or $ \chi: F^* 
\longrightarrow {\Bbb C}$ if $F$
is local) be a character of finite order cutting out a finite cyclic 
extension $L$ of $F$. Then the character 
$\chi \circ Nm: \xymatrix{
{\Bbb A}^*_E/E^* \ar[r]^{Nm} & {\Bbb A}^*_F/F^* \ar[r]^\chi & {\Bbb C}^*}$
defines the cyclic 
extension $LE$ of $E$.
\end{lemma}

\begin{cor}\label{cor6.2}
If $E$ is a quadratic extension of $F$, and $\omega: {\Bbb A}^*_E/E^* 
\longrightarrow {\Bbb C}^*$ a quadratic character defining an extension
$M$ of $E$, then
\begin{enumerate}
\item
$M$ is biquadratic over $F$ if and only if $\omega$ restricted to 
${\Bbb A}^*_F/F^*$ is trivial.
\item
$M$ is cyclic quartic over $F$ if and only if $\omega$ restricted
to ${\Bbb A}^*_F/F^*$ is $\omega_{E/F}$ where $\omega_{E/F}$
is the quadratic character on ${\Bbb A}^*_F/F^*$ defined by the
quadratic extension $E$ of $F$.
\item $M$ is non-Galois over $F$ if and only if $\omega/\omega^\sigma 
\not = 1$, and this is so if and only if $\omega$ restricted to $F^*$
is not $1$ or $\omega_{E/F}$; the restriction of $\omega$ to $F^*$ 
defines a quadratic extension, say $L^\prime$  of $F$ such that $EL^\prime$ 
is the quadratic extension of $E$ defined by $\omega\omega^\sigma$.   
\end{enumerate}
\end{cor}

\begin{proof}[Proof of corollary] One only needs to observe that $M$ is Galois
over $F$ if and only if $\omega$ is invariant under $Gal(E/F)$.
\end{proof}

Here is the main theorem regarding factorization of period 
integrals on ${\rm SL}_2$.

\begin{thm}\label{thm5.3}
Let $\pi$ be an automorphic representation of
${\rm SL}_2({\Bbb A}_E)$ contained in a cuspidal automorphic representation
$\tilde{\pi}$ of ${\rm GL}_2({\Bbb A}_E)$. Suppose that $\tilde{\pi}$
is distinguished by ${\rm GL}_2({\Bbb A}_F)$. Then
the period integral on $\pi$ is factorizable 
if $\widetilde{\pi}$ is non-monomial, or if $\widetilde{\pi}$ is monomial, and comes from
three quadratic extensions of $E$ of which exactly one is Galois over $F$. 
If $\widetilde{\pi}$ is monomial, and comes
from a unique quadratic extension, say $M$, of $E$, or comes from
three quadratic extensions of $E$ which are all Galois over $F$, 
then the period integral is not factorizable.
\end{thm}

\begin{proof} 
We recall an identity established earlier:

\begin{flushleft}
$\displaystyle{\int_{{\rm SL}_2(F)\backslash {\rm SL}_2({\Bbb A}_F)}}\phi(g)dg$ 
\end{flushleft}
\begin{flushright}
$= \displaystyle{\sum_{\begin{array}{c} \tilde{\eta}: F^*\backslash 
{\Bbb A}^*_F \rightarrow {\Bbb C}^* \\
\tilde{\eta}^2 =1 \end{array}}} \int_{{\Bbb A}^*_F {\rm GL}_2(F)
\backslash {\rm GL}_2({\Bbb A}_F)}
\phi(g)\tilde{\eta}(\det g)dg.$
\end{flushright}

In the non-monomial case, the above sum of integrals reduces to a single 
term by Corollary \ref{cor}, hence is factorizable for 
the ${\rm GL}_2$ automorphic 
representation by multiplicity one theorem for ${\rm GL}_2$,
and hence also for ${\rm SL}_2$ automorphic subrepresentations by
Lemma \ref{lemma5.1}.

If $\tilde{\pi}$ is monomial, and comes from
three quadratic extensions of $E$ of which exactly one is Galois over $F$,
then we have an isomorphism $\tilde{\pi} \cong \tilde{\pi}
\otimes \omega$, where $\omega$ is a Gr\"ossencharacter of ${\Bbb A}^*_E$
which does not restrict trivially to ${\Bbb A}^*_F$ 
(cf. Corollary \ref{cor6.2}).
Also in this case, $\tilde{\pi}$ is distinguished for exactly
two Gr\"ossencharacters of ${\Bbb A}^*_F$, namely $1$ and 
$\chi=\omega|_{_{{\Bbb A}^*_F}}$ (cf. Proposition \ref{prop}).
An isomorphism of 
${\rm GL}_2({\Bbb A}_E)$-modules 
between $\tilde{\pi}$ and $\tilde{\pi} \otimes \omega$ 
can be interpreted as an isomorphism, say $T$,
between $\tilde{\pi}$ and itself such that 
$T(gv) = \omega(\det g)gT(v)$ for all $g \in {\rm GL}_2({\Bbb A}_E)$, and 
$v \in \tilde{\pi}$. Upon modifying $T$ by a scalar, we can assume that $T$
has order 2, and splits $\tilde{\pi}$ into $\tilde{\pi}^+ 
\oplus \tilde{\pi}^-$ on which $T$
operates by $+1$  and $-1$ respectively. Since the period integral is the 
unique abstract ${\rm GL}_2({\Bbb A}_F)$-invariant linear form up to scalar, 
the $\chi$-period integral is the composite
of the period integral with $T$. The key fact is that $\pi$, being
an irreducible representation of ${\rm SL}_2({\Bbb A}_E)$, 
sits either inside $\tilde{\pi}^+$ or inside $\tilde{\pi}^-$. 
Therefore, the $\chi$-period integral on $\pi$ is a scalar multiple of the
period integral restricted to $\pi$. Hence the ${\rm SL}_2$-period integral
is factorizable by the above identity and Lemma \ref{lemma5.1}.  

In the other cases, the period integral is a 
sum of more than one linear form by Corollary \ref{cor}, 
each of which is factorizable. We argue below using Lemma \ref{lemma5.2} that 
the sum is not factorizable.  

Suppose that $\tilde{\pi}$ is `induced' from a Gr\"ossencharacter  of a 
quadratic extension $M$ of $E$. Notice that because of Lemma \ref{lemma6.5}
we can assume that $M$ is
Galois over $F$ with Galois group $({\Bbb Z}/2)^2$. However, because of
$\tilde{\pi}^\sigma \cong \tilde{\pi}^\vee$, 
if $\tilde{\pi}$ arises from a quadratic extension
$M$ of $E$, it also arises from $M^\sigma$. It is well known that 
$\tilde{\pi}$ arises either 
from 1 quadratic extension, or 3 quadratic extensions. 
It follows that one of these quadratic extensions, say $M$, of $E$ 
from which $\tilde{\pi}$ arises must be Galois over $F$. In this case,
the Gr\"ossencharacter  $\omega_{M/E}$ of ${\Bbb A}^*_E$ is $\sigma$
invariant, hence there is a character $\chi$ of ${\Bbb A}^*_E/E^*$
such that $\chi \chi^\sigma = \omega_{M/E}$. 
Clearly $\chi$ restricted to ${\Bbb A}^*_F$ 
cannot be $1$ or $\omega_{E/F}$. Since $\tilde{\pi}$ is distinguished,
$$\tilde{\pi}^\sigma \cong \tilde{\pi}^\vee \cong \tilde{\pi}^\vee 
\otimes \omega_{M/E}=\tilde{\pi}^\vee \otimes (\chi \chi^\sigma).$$
It follows that $(\tilde{\pi} \otimes \chi^{-1})^\sigma \cong 
(\tilde{\pi} \otimes \chi^{-1})^\vee$, 
and therefore $\tilde{\pi}$ is $\chi$ or $\chi\omega_{E/F}$-distinguished,
which after perhaps changing the choice of $\chi$ with 
$\chi \chi^\sigma = \omega_{M/E}$, we can assume that $\tilde{\pi}$ 
is $\chi$ distinguished (besides being distinguished).
We note that this implies, in particular, that $\chi$ restricted to 
${\Bbb A}^*_F$ is of order $2$ (we have already noted earlier that $\chi$
restricted to ${\Bbb A}^*_F$ is not trivial).

Observe now that the quadratic character
$\chi$ restricted to ${\Bbb A}^*_F$ 
defines a quadratic extension $L$ of $F$, and from the equality
$\chi \chi^\sigma = \omega_{M/E}$, $M$ is a biquadratic extension,
$M = LE$ of $F$. Assume that $\tau$ is the
non-trivial automorphism of $M$ over $E$, and abusing notation,
let $\sigma$ be the non-trivial automorphism of $M$ over $L$.

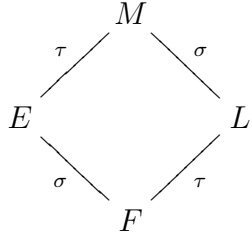
$$
\xymatrix{
& M \ar@{-}[dl]_\tau \ar@{-}[dr]^\sigma & \\
E \ar@{-}[dr]_\sigma & & L \ar@{-}[dl]^\tau \\
& F& 
}
$$
Suppose that $\tilde{\pi}$ arises from a Gr\"ossencharacter 
$\mu$ of ${\Bbb A}_M^*/M^*$,
and is distinguished. Therefore,
$$\tilde{\pi}^\sigma \cong \tilde{\pi}^\vee,$$
which assuming $\tilde{\pi}$ arises from $\mu$ implies that 
$${\rm Ind}_{W_M}^{W_E} (\mu^\sigma) = {\rm Ind}_{W_M}^{W_E}(\mu^{-1}).$$
This implies that either,
\begin{eqnarray}
\mu^\sigma & = & \mu^{-1},  
\end{eqnarray}
or, 
\begin{eqnarray}
\mu^{\sigma \tau}&  = & \mu^{-1}.
\end{eqnarray}

Defining $L_1$ to be the field fixed by $\sigma\tau$, we 
note that in case $(3)$, the character $\mu /\mu^\tau$ is trivial on 
 ${\Bbb A}_{L_1}^*$, 
because in this case $\mu \mu^{\sigma \tau}=1$,
therefore $\mu$ restricted to ${\Bbb A}_{L_1}^*$ is either $1$
 or $\omega_{M/L_1}$. Therefore $\mu$ and $\mu^\tau$ have the same restriction
to ${\Bbb A}_{L_1}^*$, proving our claim. 

In case $(2)$, we claim that $\mu /\mu^\tau$ restricted to  
 ${\Bbb A}_{L_1}^*$ cannot be $\omega_{M/L_1}$. If
$\mu /\mu^\tau$ restricted to 
 ${\Bbb A}_{L_1}^*$ was $\omega_{M/L_1}$, then in particular
$(\mu /\mu^\tau) (\mu /\mu^\tau)^{\sigma \tau} =1$. Since we are
in case $(2)$, $\mu^\sigma = \mu^{-1}$. Therefore, the
condition $(\mu /\mu^\tau) (\mu /\mu^\tau)^{\sigma \tau} =1$
becomes $(\mu /\mu^\tau)^2 = 1$. If $\mu /\mu^\tau$ restricted to  
 ${\Bbb A}_{L_1}^*$ is $\omega_{M/L_1}$, then by Corollary 
6.7, the quadratic extension of $M$ 
defined by $\mu/\mu^\tau$, call it $M_1$,
 is a cyclic quartic extension of $L_1$. 
Now we note the following elementary lemma whose proof is omitted.
\begin{lemma}\label{lemma_el}
Let $N$ be an abelian normal subgroup of a group $G$ with $G/N$ cyclic.
Assume that the action of $G/N$ on $N$ via inner conjugation is trivial.
Then $G$ is abelian.
\end{lemma}

We apply the above lemma to $G=Gal(M_1/F)$ which contains 
$N=Gal(M_1/E)=({\Bbb Z}/2)^2$ as an abelian normal subgroup on which
$Gal(E/F)$ acts trivially since we are in the situation in which
all the three quadratic extensions of $E$ from which $\tilde{\pi}$ arises
are Galois over $F$. Thus we find that $G=Gal(M_1/F)$ is abelian.
Because of Lemma \ref{lemma6.5}, $G$ does not contain ${\Bbb Z}/4$ as
a quotient and hence neither as a subgroup.
This implies that the Galois
group of $M_1$ over $L_1$ cannot be ${\Bbb Z}/4$, proving our
claim that $\mu /\mu^\tau$ restricted to  
 ${\Bbb A}_{L_1}^*$ cannot be $\omega_{M/L_1}$.

Before proceeding further, we note the following lemma which is 
at the basis of our proof of non-factorization of period integral. 
This is part of 
case 3 of Theorem 1.3 of our paper \cite{anand}. It can be easily
proved by a direct analysis of the ${\rm GL}_2(F_v)$ action on 
${\Bbb P}^1(E_v)$.
\begin{lemma}\label{lemma6.6}
Let $E_v$ be a quadratic extension of a local field $F_v$,
and $\pi = Ps(\chi_1,\chi_2)$ a principal series representation
of ${\rm GL}_2(E_v)$. Then if $\chi_1=\chi_2$, $\pi$ remains irreducible
when restricted to ${\rm SL}_2(E_v)$, and is $\nu$ distinguished for two
characters $\nu$ of $F_v^*$.
\end{lemma}

In what follows, we will be doing some local analysis for which
we assume that all
our places in consideration in $L_1$ or $M$ are unramified
over the corresponding place in $F$, and the character $\mu$ 
is unramified at these places.

 We note that there are infinitely many primes in $L_1$ which are inert
in $M$. The prime in $F$ below such a 
prime in $L_1$ has the property that it is 
inert in both $L$ and $E$, and split in $L_1$.
We abuse notation to denote the pair of places in $M$ as well
as in $L_1$ as $(v_1,v_2)$. Since we are going to use only unramified
characters, this should not cause any confusion.

 If the local components of $\mu$ 
at $(v_1,v_2)$ is $(\mu_1,\mu_2)$,  $\mu/\mu^\tau$ looks like 
$(\mu_1/\mu_2,\mu_2/\mu_1)$ at this pair of places. 

In case $(3)$, since the character $\mu /\mu^\tau$ is trivial on 
 ${\Bbb A}_{L_1}^*$, in particular the pair of characters
$(\mu_1/\mu_2,\mu_2/\mu_1)$ is trivial, hence $\mu_1 = \mu_2$. 

In case $(2)$, we know that $\mu /\mu^\tau$ restricted to  
 ${\Bbb A}_{L_1}^*$ is a certain quadratic character which is
not $\omega_{M/L_1}$. 
Either, $\mu /\mu^\tau$ restricted to  
 ${\Bbb A}_{L_1}^*$ is the trivial character, in which case
places of $L_1$ which are inert in $M$ automatically give
$\mu_1 = \mu_2$, or the quadratic extension of $L_1$
defined by  $\mu /\mu^\tau$ restricted to  
 ${\Bbb A}_{L_1}^*$ is distinct from $M$, and together with $M$ 
gives a Galois
extension of $F$ with Galois group $({\Bbb Z}/2)^3$.
Applying Cebotaraev density
theorem, we once again find that there are infinitely many
primes of $L_1$ which are inert in $M$ where the restriction of
$\mu /\mu^\tau$ is trivial. Hence, once again $\mu_1 = \mu_2$.

Thus Lemma \ref{lemma6.6} applies, and which in conjunction 
with Lemma \ref{lemma5.2}
implies that the period integral is not factorizable, 
completing the proof of Theorem \ref{thm5.3}.
\end{proof}

\noindent
\begin{remark} 
Observe that the above theorem can be viewed as an analogue
of Theorem 1.2 of \cite{anand}. The cases where the period integral is Eulerian
are exactly the global analogues of the 
cases in Theorem 1.2 of \cite{anand} where the 
space of local invariant forms has multiplicity
one. Note that this analogy holds in the context of
Jacquet's result too \cite{jacquet3}.
There the symmetric space $({\rm Res}_{E/F}{\rm GL}(3),{\rm U}(3))$ 
has the property that locally over a $p$-adic field, the 
space of ${\rm U}(3)$-invariant linear forms on a 
supercuspidal representation of ${\rm GL}_3(E)$ has multiplicity at most one, 
and over global fields, the period integral is factorizable for
cuspidal representations.
\end{remark}

\begin{remark}
For a reductive algebraic group $G$ over a local field $k$, 
$K$ a separable quadratic extension of $k$, and $\pi$ an irreducible admissible
representation of $G(K)$, it makes sense to study the dimension of the
space of $G(k)$-invariant forms $\ell:\pi \rightarrow {\Bbb C}$.
It is reasonable to expect that this dimension is always finite. In the global
study, since the linear form is fixed to be the period integral, there is
no obvious global analogue of the concept of the dimension
of $G(k)$-invariant forms. However Theorem \ref{thm5.3}, seen in the light
of the corresponding local result,
Theorem 1.2 of \cite{anand}, suggests
a reasonable global analogue to be the smallest positive integer
$d$ such that the period integral can be written as a sum of $d$ factorizable
linear forms. With this notion, we can go a step further in 
Theorem \ref{thm5.3} to say that in the cases in which the period integral is 
not factorizable, it is a sum of two or four factorizable linear forms 
depending on whether the representation comes from a unique quadratic extension
of $E$ or three quadratic extensions of $E$ which are all Galois over $F$.
We omit the details of this calculation. It is curious to note that
not only is $d$ finite for $SL_2$, it has a very similar structure
to the dimension of the space of local invariant forms. Understanding
these local and global dimensions in general seems a very interesting
problem. In this connection, we mention the work of  Lapid and Rogawski 
\cite{lapid3}
which computes the a period of an Eisenstein series on $GL_3$ as a sum of
factorizable functionals, its recent generalization by
Omer Offen \cite{offen}, as well as the 
earlier work of 
Jacquet \cite {jacquet3} for $GL_3$, and its recent 
generalization to $GL_n$.

\end{remark}

\section{Pseudo-Distinguishedness} 
If an automorphic  representation $\pi = \otimes\pi_v$ 
of  ${\rm GL}_2({\Bbb A}_E)$ 
has the property that $\pi_v$ is  distinguished by ${\rm SL}_2(F_v)$ 
at all places $v$ of $E$, then there are characters $\chi_v$ 
of $E_v^*$ such that 
$$(\pi_v \otimes \chi_v)^\sigma \cong (\pi_v \otimes \chi_v)^\vee.$$
Thus at all places $v$ of $E$, $\pi_v^\sigma$ and 
$\pi_v^\vee$ differ by a character of $E_v^*$. By the multiplicity
one theorem of Ramakrishnan, cf. \cite{dinakar1}, this implies that 
$$\pi^\sigma \cong \pi ^\vee \otimes \chi$$
for a character $\chi$ of ${\Bbb A}^*_E/E^*$.

The aim of this section is to classify representations
$\pi$ of ${\rm GL}_2({\Bbb A}_E)$ such that 
\begin{eqnarray} \label{eq1}
\pi^\sigma \cong \pi ^\vee \otimes \chi,
\end{eqnarray}
for a character $\chi$ of ${\Bbb A}^*_E/E^*$ which we assume fixed in 
this section, 
and which is not Galois invariant.
We call such representations {\it pseudo-distinguished}. Although,
we write the arguments below for $\pi$ an irreducible admissible
representation of ${\rm GL}_2(E)$, $E$ a local field, exactly
the same argument works in the case of automorphic forms over 
global fields. We note that Lapid and Rogawski, cf. \cite{lapid2}, 
have also done an analogous study,
of classifying $\pi$ with $\pi^\sigma \cong \pi \otimes \chi$,   
via the methods of trace formula. 

We note that if $\chi$ is Galois invariant, then we can write
$\chi$ as $\chi = \alpha \cdot \alpha^\sigma$ for a character
$\alpha $ of $E^*$, and therefore equation (\ref{eq1}) reduces 
after twisting by a character to $\pi^\sigma \cong \pi^\vee$,
studied in the theory of distinguished representations.

By applying $\sigma$ to (\ref{eq1}), and rewriting, we find,
\begin{eqnarray}\label{eq2}
\pi^\sigma \cong \pi ^\vee \otimes \chi^\sigma. 
\end{eqnarray}
Therefore from (\ref{eq1}) and (\ref{eq2}), if $\chi^\sigma \not = \chi$,
then $\chi^\sigma/\chi$ is a quadratic character, say $\omega$,
of $E^*$, and $\pi$ has a self-twist by $\omega$, implying 
that $\pi$ is a monomial representation arising from 
a character $\mu$ of the quadratic extension
$M$ of $E$ defined by by $\omega$: $\pi = {\rm Ind}_{W_M}^{W_E} \mu$,
where $W_M$ and $W_E$ are respectively the Weil groups of $M$ and $E$.
Since $\chi(x/x^\sigma) = \omega(x)$, $\omega$ restricted to ${\Bbb A}^*_F$ 
is trivial. Therefore by Corollary \ref{cor6.2}, 
$M$ is a biquadratic extension of $F$, say $M=EL$ with $L$ a quadratic
extension of $F$. Assume that $\tau$ is the
non-trivial automorphism of $M$ over $E$, and abusing notation,
let $\sigma$ be the non-trivial automorphism of $M$ over $L$.

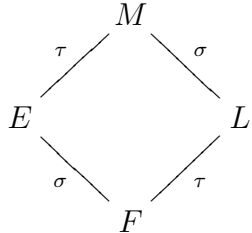
$$
\xymatrix{
& M \ar@{-}[dl]_\tau \ar@{-}[dr]^\sigma & \\
E \ar@{-}[dr]_\sigma & & L \ar@{-}[dl]^\tau \\
& F& 
}
$$

Once again, condition $(4)$,
$${\rm Ind}_{W_M}^{W_E} (\mu^\sigma) = {\rm Ind}_{W_M}^{W_E}(\mu^{-1}) 
\otimes \chi,$$
implies that either,
\begin{eqnarray}
\mu^\sigma & = & \mu^{-1} \chi \chi^\tau  
\end{eqnarray}
or, 
\begin{eqnarray}
\mu^{\sigma \tau}&  = & \mu^{-1} \chi \chi^\tau .
\end{eqnarray}
Let us consider the first case, as the other case is similar. In this case,
\begin{eqnarray}
\mu \mu^\sigma = \chi \chi^\tau.
\end{eqnarray}

Note that since $\chi/\chi^\sigma$ is of order 2 defining $M$ over
$E$, the character $\chi \chi^\tau$ of ${\Bbb A}^*_M/M^*$ 
is $\sigma$-invariant. Therefore there exists a character $\tilde{\mu}$
of ${\Bbb A}^*_L/L^*$ such that 
$$\tilde{\mu}(Nm x) = (\chi \chi^\tau)(x),$$
for all $x \in {\Bbb A}^*_M/M^*$ where the norm is taken from 
${\Bbb A}^*_M/ M^*$ to $ {\Bbb A}^*_L/L^*$. 
Any extension $\mu$ of $\tilde{\mu}$ from 
${\Bbb A}^*_L/ L^*$ to
$ {\Bbb A}^*_M/M^*$ satisfies condition $(6)$.  
Further, it can be seen that $\tilde{\mu}$
is unique up to the automorphism of $L$ over $F$ 
(since the automorphism takes $\tilde{\mu}$ to $\tilde{\mu}\omega_{M/E}$), 
proving the following proposition.
\begin{prop}\label{prop6.3}
Given a character $\chi$ of ${\Bbb A}^*_E/ E^*$ such that $\chi/\chi^\sigma$
is a character of order 2 of ${\Bbb A}^*_E/ E^*$, defining a 
quadratic extension
$M$ of $E$ which is of the form $M= LE$ 
where $L$ is a quadratic extension of $F$, there
exists a character $\tilde{\mu}$ of ${\Bbb A}^*_L/ L^*$ such that the
characters $\chi$ and $\tilde{\mu}$ restricted to 
${\Bbb A}^*_M/ M^*$ (via the norm maps to ${\Bbb A}^*_E/ E^*$ and 
${\Bbb A}^*_L/ L^*$) are
the same, i.e., $\chi(x\cdot x^\tau) = \tilde{\mu}(x \cdot
x^\sigma)$ for all $x \in {\Bbb A}^*_M/  M^*$. Such a character $\tilde{\mu}$
is unique up to the Galois automorphism of $L$ over $F$
(since the automorphism takes $\tilde{\mu}$ to $\tilde{\mu}\omega_{M/E}$), 
and any extension $\mu$ of $\tilde{\mu}$ from ${\Bbb A}^*_L/ L^*$ to $
 {\Bbb A}^*_M/M^*$
gives rise to an isomorphism
$$\pi^\sigma \cong \pi ^\vee \otimes \chi, $$
for $\pi = {\rm Ind}_{W_M}^{W_E} (\mu)$.
\end{prop}

\begin{remark} 
Since characters of ${\Bbb A}^*_M/ M^*$ extending
a given character of ${\Bbb A}^*_L/ L^*$ are --after fixing one such 
character-- 
in bijective correspondence with characters of $ {\Bbb A}^*_M/M^*$ trivial on
$ {\Bbb A}^*_L/L^*$, one can state the proposition in a more suggestive way
as follows: representations $\pi$ of ${\rm GL}_2(E)$ with 
$\pi^\sigma \cong \pi ^\vee \otimes \chi, $ with $\chi/\chi^\sigma$
cutting out a quadratic extension $M=LE$ of $E$ are in bijective
correspondence with representations of
${\Bbb A}^*_M/ M^*$ distinguished by ${\Bbb A}^*_L/ L^*$.
\end{remark}

\section{Locally but not globally distinguished II}
In this section, 
we construct an automorphic representation
$\pi=\otimes \pi_v$  of ${\rm SL}_2({\Bbb A}_E)$ which is
 abstractly ${\rm SL}_2({\Bbb A}_F)$ distinguished
but none of the elements in the global $L$-packet 
determined by $\pi$ is  distinguished by ${\rm SL}_2({\Bbb A}_F)$.
We achieve this by the following steps.
\begin{enumerate}
\item 
We construct a pseudo-distinguished representation $\tilde{\pi} 
= \otimes \tilde{\pi}_v$ of ${\rm GL}_2({\Bbb A}_E)$ with
$$\tilde{\pi}^\sigma \cong \tilde{\pi}^\vee \otimes \chi,$$
for a Gr\"ossencharacter $\chi$ with $\chi^\sigma \not = \chi$.
Our representation $\tilde{\pi}$ will be monomial arising from 
exactly one quadratic extension of $E$, and hence there is exactly
one non-trivial quadratic character $\omega$ such that 
$$ \tilde{\pi} \cong \tilde{\pi} \otimes \omega.$$
This implies that the only  Gr\"ossencharacters $\alpha$ with
$$ \tilde{\pi}^\sigma \cong \tilde{\pi}^\vee \otimes \alpha,$$
are $\chi$ and $\chi^\sigma$, and in particular, 
there are none with  $\alpha^\sigma = \alpha$.
\item 
Step $(1)$ implies that such an automorphic representation
of ${\rm GL}_2({\Bbb A}_E)$ is not $\nu$-distinguished with respect to
${\rm GL}_2({\Bbb A_F})$ for any Gr\"ossencharacter $\nu$ of ${\Bbb A}_F^*$,
and hence by Proposition \ref{prop3.3}, none of the members of the $L$-packet
of automorphic forms determined by $\tilde{\pi}$ is ${\rm SL}_2({\Bbb A}_F) 
$-distinguished.
\item 
We next ensure that $\tilde{\pi}$ is locally distinguished
(with respect to some character of $F_v^*$) 
at all the places $v$
of $F$, and hence is abstractly distinguished by $SL_2({\Bbb A}_F)$.
\item 
At a place $v$ of $F$ which splits as $(v_1,v_2)$ in $E$,
the condition 
$$\tilde{\pi}^\sigma \cong \tilde{\pi}^\vee \otimes \chi,$$ 
amounts to   $(\tilde{\pi}_{v_1},\tilde{\pi}_{v_2})^\sigma
\cong (\tilde{\pi}_{v_1},\tilde{\pi}_{v_2})^\vee 
\otimes (\chi_1,\chi_2)$. This
is equivalent to 
\begin{eqnarray*}
\tilde{\pi}_{v_2} & \cong & \tilde{\pi}_{v_1}^\vee\cdot \chi_1 \\
\tilde{\pi}_{v_1} & \cong & \tilde{\pi}_{v_2}^\vee \cdot \chi_2. 
\end{eqnarray*}
This implies in particular that $\tilde{\pi}_{v_1} \otimes \tilde{\pi}_{v_2}$
has $\chi_2$-invariant linear form. Therefore in the $L$-packets
determined by $\tilde{\pi}_{v_1}$ and $\tilde{\pi}_{v_2}$,
there are representations ${\pi}_{v_1}$ and ${\pi}_{v_2}$ of
${\rm SL}_2(F_v)$ such that $\pi_{v_1} \cong \pi^\vee_{v_2}$. Thus for places
of $F$ which are split in $E$, local distinguishedness is automatic
from the pseudo-distinguishedness condition 
$\tilde{\pi}^\sigma \cong \tilde{\pi}^\vee \otimes \chi.$
\item 
For a place $v$ of $E$ which is inert over $F$, we will 
ensure that the local representation $\tilde{\pi_v}$ is either 
unramified, or comes from an unramified character of a quadratic
extension, say $M_v$ of $E_v$ which is Galois over $F_v$. In the latter case,
$\tilde{\pi_v}$ is the
 principal series representation of the form $Ps(\chi,\chi \omega)$
where $\chi$ is an unramified character of $E_v$, and $\omega$ is the
quadratic character defining the quadratic extension $M_v$ 
of $E_v$, and is 
therefore invariant under the automorphism of $E_v$ over $F_v$.
From Lemma \ref{lemma8.1} below,  $\tilde{\pi_v}$
is ${\rm SL}_2(F_v)$-distinguished.
\item 
At places $v$ of $E$ at which $\tilde{\pi}_v$ is unramified, and
$v$ itself is unramified over $F$, all the members of the $L$-packet
determined by $\tilde{\pi}_v$ are distinguished by ${\rm SL}_2(F_v)$.
This easily follows as under these conditions 
${\rm GL}_2(F_v)$ operates transitively on the
$L$-packet of ${\rm SL}_2(E_v)$ determined by $\tilde{\pi}_v$.
\item 
By steps (4),(5),(6), there are ${\rm SL}_2({\Bbb A}_E)$ components
of $\tilde{\pi} = \otimes \tilde{\pi}_v$ which are ${\rm SL}_2({\Bbb A}_F)
$-distinguished. Because of the flexibility offered by step 6, we can
assume that these are even automorphic.
\end{enumerate}

The proof of the following elementary lemma follows from
Proposition 2.3 of [A-P].

\begin{lemma}\label{lemma8.1}
For a separable quadratic extension $K$ of a non-archimedean 
local field $k$ with the nontrivial Galois automorphism $\sigma$ 
of $K$ over $k$, a principal series representation $Ps(\chi_1,\chi_2)$ of
$GL_2(K)$ is distinguished by $SL_2(k)$ if and only if either 
$(\chi_1\chi^{-1}_2)|_{k^*} = 1$, or $(\chi_1\chi^{-1}_2)
= (\chi_1\chi^{-1}_2)^\sigma$.
\end{lemma}

Here is the main theorem of this section.

\begin{thm}\label{thm7.1}
There exists a cuspidal automorphic representation $\pi$ of
${\rm SL}_2({\Bbb A}_{E})$ for $E= {\Bbb Q}(\sqrt {-1})$  
which is 
locally distinguished with respect to $SL_2({\Bbb A}_{\Bbb Q})$ 
at all the places
of ${\Bbb Q}$, but for which none of the members 
of its $L$-packet  is globally distinguished.
\end{thm}

\begin{proof} We will construct a cuspidal
representation $\tilde{\pi}$ of ${\rm GL}_2({\Bbb A}_{E})$
which is pseudo-distinguished  by a character $\chi$ 
of ${\Bbb A}_{E}^*/{E}^*$ with $\chi \not = \chi^\sigma$, and  unramified
at all places of ${E}$. Such a representation $\tilde{\pi}$ 
is locally ${\rm SL}_2$-distinguished at all places of ${\Bbb Q}$. 
The representation
$\tilde{\pi}$ will be a monomial representation, coming from exactly one
quadratic extension of $E$. As we have seen above (in step 2),
the $L$-packet of automorphic representations of $ {\rm SL}_2({\Bbb A}_{E})$
determined by $\tilde{\pi}$ has no globally 
${\rm SL}_2({\Bbb A}_{\Bbb Q})$-distinguished member. We now construct
the specific example.

Let $L= {\Bbb Q}(\sqrt{-257})$, and $L'=  {\Bbb Q}(\sqrt{257})$. 
From the tables in [B-S], the class group $C_L$ of $L$ is ${\Bbb Z}/16$,
and the class group $C_{L'}$ of $L'$ is ${\Bbb Z}/3$. Let $M = 
 {\Bbb Q}(\sqrt{257}, \sqrt{-257})$ be the unique quadratic unramified 
extension of $L$. The natural map from the class group
$C_L$ to the class group
$C_M$ has ${\Bbb Z}/2$ as its kernel,  and ${\Bbb Z}/8$ as its image.
Further, the image of the 
natural map from $C_{L'}$ to $C_M$ is ${\Bbb Z}/3$. 

Let $\tau$ be the automorphism 
of $M$ over ${\Bbb Q}$ which is nontrivial on both $L$ and $L'$, and we abuse
notation to denote its restriction to $L$ or $L'$ also by $\tau$. 
Further we let $\sigma$ be the nontrivial automorphism of $M$ over ${\Bbb Q}$
which is trivial on $L$.

We will be constructing an unramified character 
$\mu^\prime$ of  ${\Bbb A}_L^*/L^*$, which is the same as a character of
$C_L$,
such that $\mu^\prime/\mu^{\prime\tau}$ is a quadratic character, and hence
defines the quadratic unramified 
extension $M$ of  $L$. For
any extension $\mu$ of $\mu^\prime$ to ${\Bbb A}_M^*/M^*$, it is easy
to see that $\mu \mu^\sigma$ is a $\tau$-invariant character
on  ${\Bbb A}_M^*/M^*$, and hence there is a character $\tilde{\chi}$
of  ${\Bbb A}_M^*/M^*$ such that
$$\mu \mu^\sigma = \tilde{\chi} \tilde{\chi}^\tau.$$

Denote
the restriction of $\tilde{\chi}$ to ${\Bbb A}_E^*/E^*$ by $\chi$.
It can be checked that $\chi/\chi^\sigma$ is the quadratic character
of $E$ defining the quadratic extension $M$ of $E$.

It is
easy to see that the action of $\tau$ on $C_L$ 
and also on $C_{L'}$ is $x \rightarrow -x$, and hence also on the 
image of these groups in $C_M$.

Let $\mu'$
be a character of $C_L$  of order 4. Such a character
 $\mu'$ is trivial on the kernel of the map from
$C_L$ to $C_M$. Since $\tau$ acts by $ x \rightarrow -x$ on $C_L$ and
$\mu'$ is of order 4,
$\mu'/\mu'^\tau$ is a character of order 2, hence defines the unique
quadratic unramified extension of $L$ which we are denoting by $M$.
The character $\mu'$ being trivial on the kernel of the map 
from $C_L$ to $C_M$ extends to a character $\mu$ of $C_M$ which we take to be
nontrivial on the ${\Bbb Z}/3$ coming from $L'=  {\Bbb Q}(\sqrt{257})$.
Therefore $\mu/\mu^\tau$ is not of order 2. 
Let $\tilde{\pi}$ be 
the cuspidal representation  
on ${\rm GL}_2({\Bbb A}_E)$ obtained by `inducing' the character $\mu$ of 
 ${\Bbb A}_M^*/M^*$, which by the condition
$\mu \mu^\sigma = \tilde{\chi} \tilde{\chi}^\tau,$ will be
pseudo-distinguished for the character $\chi$. The 
representation $\pi$ sought after in the statement of the theorem is 
obtained from $\tilde{\pi}$ as in the paragraph before the theorem.

\end{proof}


\end{document}